\documentclass[12pt]{article}
\usepackage{amssymb,amsmath}
\def\hang{\hangindent\parindent}
\def\tex#1{\indent\llap{[#1]\enspace}\ignorespaces}
\def\re{\par\hang\tex}
\def\F{{\mathbb{F}}}

\def\D{\Delta}
\def\det{{\rm det}}

\def\LL{{\mathcal L}}

\def\sc{\scriptstyle}
\def\ssc{\scriptscriptstyle}
\def\dis{\displaystyle}
\def\cl{\centerline}

\def\vs{\vspace*}

\def\ni{\noindent}

\def\Z{\mathbb{Z}{\ssc\,}}

\def\C{\mathbb{C}{\ssc\,}}

\numberwithin{equation}{section}
\newtheorem{theo}{Theorem}[section]

\newtheorem{lemm}[theo]{Lemma}

\begin{document}
\cl{{\large \bf Classification of $\Z^2$-graded modules of the
intermediate series}}\cl{{\large\bf over a Lie algebra of Block
type} \footnote {Supported by NSF grant no.~10471091 of China, the
grant ``One Hundred Talents Program'' from University of Science and
Technology of China.}} \vs{6pt}

\cl{Xiaoqing Yue$^{\dag}$,  \ \ Yucai Su$^\ddag$} \cl{\footnotesize
$^\dag$Department of Mathematics, Tongji University, Shanghai
200092, China} \cl{\footnotesize$^\ddag$Department of Mathematics,
University of Science and Technology of China, Hefei 230026, China}
\cl{\footnotesize E-mail: yuexqtj@21cn.com, ycsu@ustc.edu.cn}
\vs{6pt} \vs{10pt}\par

{\small \parskip .005 truein \baselineskip 3pt \lineskip 3pt

\noindent{{\bf Abstract.} Let $\LL$ be the Lie algebra of Block type
over a field $\F$ of characteristic zero, defined with basis
$\{L_{i,j}\,|\,i,j\in\Z \}$ and relations
$[L_{i,j},L_{k,l}]=((j+1)k-(l+1)i)L_{i+k,j+l}.$ Then $\LL$ is
$\Z\times\Z$-graded.  In this paper, $\Z^2$-graded $\LL$-modules of
the intermediate series are classified.
 \vs{5pt}

}

}
\parskip .001 truein\baselineskip 8pt \lineskip 8pt

\vs{6pt}
\par

\cl{\bf\S1. \
Introduction}\setcounter{section}{1}\setcounter{equation}{0} Let
$\LL$ be the {\it Lie algebra of Block type} defined over a field
$\F$ of characteristic zero, with basis $\{L_{i,j}\,|\, i,j\in\Z\}$
and \vs{-4pt}relations
\begin{equation}\label{def1}
[L_{i,j},L_{k,l}]=((j+1)k-(l+1)i)L_{i+k,j+l}\mbox{ \ \ for \
}i,j,k,l\in\Z\vs{-4pt}.
\end{equation}
This Lie algebra is closely related to the Virasoro algebra,
Virasoro-like algebra. It is also a special case of Cartan type $S$
Lie algebra or  Cartan type $H$ Lie algebra. Partially due to these
facts, the study of this Lie algebra or its analogs has recently
attracted some authors' attentions (see e.g., [DZ, LT, S1, S2, WS,
X1, X2, ZZ, ZM]). It is well known that although Cartan type Lie
algebras have a long history, their representation theory is however
far from being well developed.

In order to better understand the representation theory of Cartan
type Lie algebras, it is very natural to first study representations
of special cases of Cartan type Lie algebras. The author in [S1, S2]
presented a classification of the so-called {\it quasifinite
modules} (which are simply $\Z$-graded modules with finite
dimensional homogenous subspaces) over some Block type Lie algebras.
In particular, it is shown that quasifinite modules with dimensions
of homogenous subspaces being uniformly bounded are all
trivial---this is due to the crucial fact that all $\Z$-graded
homogenous subspaces of these Block type Lie algebras are themselves
infinite dimensional. On the other hand, the authors in [WS]
considered Verma type modules of the above Lie algebra $\LL$.
However it turns out that these Verma type modules, regarded as
$\Z$-graded modules, are all with infinite dimensional homogenous
subspaces (except the subspace spanned by the generator which has
dimension $1$). The author in [Z] gave a classification of weight
modules with $1$-dimensional weight spaces over (generalized) Witt
algebras (the first family of Cartan type Lie algebras). Since our
focus in the present paper is the Block type Lie algebras, in above
we only listed very few results on representations of Cartan type
Lie algebras (and of course there are many other results that we did
not mention), but one may have already realized that representations
of Cartan type Lie algebras are not as easy as one might expect.

The aim in the present paper is to give a classification of
$\Z^2$-graded $\LL$-modules with $1$-dimensional homogenous
subspaces (which are usually called {\it modules of the intermediate
series}). Our next goal is to give a classification of all
$\Z^2$-graded $\LL$-modules with finite dimensional homogenous
subspaces in some due time. The main techniques used in this paper
are those developed in [S3, S4] for representations of higher rank
Virasoro algebras. It seems to us that the techniques developed in
[S3, S4] are useful in dealing with problems of representations of
Lie algebras related to the Virasoro algebra.

The main result of this paper is the following.
\begin{theo}\label{theo}
Suppose $V=\oplus_{i,j\in \Z}V_{i,j}$ is an indecomposable $\LL$-module of the
intermediate series $($i.e., $\dim\,V_{ij}\le 1$ for all
$i,j\in\Z)$,  then $V$ is one of  $A_{a,b},\,A,\,B,\,C$ or one of
their quotient modules,
 for some $a,b\in\F$, where $A_{a,b},\,A,\,B,\,C$ all have
 basis $\{v_{k,l}\ |\ k,l\in\Z\}$, such that
\begin{eqnarray}\label{A-a-b}
\!\!\!\!\!\!\!\!\!\!\!\!&\!\!\!\!\!\!\!\!\!\!\!\!\!\!\!& A_{a,b}:
L_{i,j}v_{k,l}=((j+1)(a+k)+(b-1-l)i)v_{i+k,j+l}\ ;
\\\!\!\!\!\!\!\!\!\!\!\!\!&\!\!\!\!\!\!\!\!\!\!\!\!\!\!\!&
A:\ \ \
L_{i,j}v_{k,l}=((j+1)k-(l+1)i)v_{i+k,j+l}+\delta_{i+k,0}\delta_{l+j,-2}iv_{0,-2}\
\mbox{ for \ }(k,l)\not=(0,-2),
\nonumber\\\!\!\!\!\!\!\!\!\!\!\!\!&\!\!\!\!\!\!\!\!\!\!\!\!\!\!\!&\
\ \ \ \ \ \ \ L_{i,j}v_{0,-2}=0\ ;\ \label{B-a-b}
\\\!\!\!\!\!\!\!\!\!\!\!\!&\!\!\!\!\!\!\!\!\!\!\!\!\!\!\!&
B:\ \ \ L_{i,j}v_{k,l}=((j+1)k-(l+1)i)v_{i+k,j+l}\
\nonumber\\\!\!\!\!\!\!\!\!\!\!\!\!&\!\!\!\!\!\!\!\!\!\!\!\!\!\!\!&\
\ \ \ \ \ \ \ \ \ \ \ \ \ \ \ \ \ \ \ \ \ \ \ \ \ \ \ \ \ \mbox{ for
\ } (k,l)\not=(0,-1),(0,-2)\ \mbox{ and \ } (i+k,j+l)\not=(0,-1),
\nonumber\\\!\!\!\!\!\!\!\!\!\!\!\!&\!\!\!\!\!\!\!\!\!\!\!\!\!\!\!&\
\ \ \ \ \ \ \ L_{i,j}v_{-i,-j-1}=0\ , \ \ \ \ \ \ \
L_{i,j}v_{0,-1}=iv_{i,j-1}\ , \ \ \ \ \ \ \
L_{i,j}v_{0,-2}=iv_{i,j-2}\ ;\label{B-a-b-c}
\\\!\!\!\!\!\!\!\!\!\!\!\!&\!\!\!\!\!\!\!\!\!\!\!\!\!\!\!&
C:\ \ \
L_{i,j}v_{k,l}=((j+1)k-(l+1)i)v_{i+k,j+l}+\delta_{i+k,0}\delta_{l+j,-2}iv_{0,-2}\
\nonumber\\\!\!\!\!\!\!\!\!\!\!\!\!&\!\!\!\!\!\!\!\!\!\!\!\!\!\!\!&\
\ \ \ \ \ \ \ \ \ \ \ \ \ \ \ \ \ \ \ \ \ \ \ \ \ \ \ \ \ \mbox{ for
\ } (k,l)\not=(0,-1),(0,-2)\ \mbox{ and \ } (i+k,j+l)\not=(0,-1),
\nonumber\\\!\!\!\!\!\!\!\!\!\!\!\!&\!\!\!\!\!\!\!\!\!\!\!\!\!\!\!&\
\ \ \ \ \ \ \ L_{i,j}v_{-i,-j-1}=0\ , \ \ \ \ \ \ \
L_{i,j}v_{0,-1}=iv_{i,j-1}\ , \ \ \ \ \ \ \
L_{i,j}v_{0,-2}=0.\label{B-a-b-d}
\end{eqnarray}
\end{theo}

The above results show that, unlike the case for the Virasoro
algebra, for the Lie algebra $\LL$, somehow it is strange that the
module
 $A_{a,b}$  only has three isolated
 deformations $A,\,B,\,C$ without a parameter (cf.~(\ref{condition})--(\ref{condition++})).

 To describe the structure of modules of the intermediate series,
we introduce the following notation: For a module $V$ and two
subsets $V_1,\,V_2$, if $V_2$ can be generated by $V_1$, namely,
$V_2\subset U(\LL)V_1$ (where $U(\LL)$ is the universal enveloping
algebra of $\LL$), then we denote $V_1\to V_2$ or $V_2\leftarrow
V_1$. Then the following result can be easily verified.
\begin{theo}\label{TH2}
\begin{itemize}\item[{\rm(i)}] $A_{a,b}$ is irreducible if and only if $a\notin\Z$
or $a\in\Z$ and $b\notin\Z$.
\item[{\rm(ii)}] If $a,b\in\Z$, then $A_{a,b}$ has 3 composition factors with the following
structure,
\begin{eqnarray}
V_2\to V_3\to V_1,
\end{eqnarray}
where $V_1=\F v_{-a,b-1},\,V_2=\F v_{-a,b-2},\,V_3= {\rm
span\,}{\{v_{k,l}\,|\,(k,l)\not=(-a,b-1),(-a,b-2)\}}.$
\item[{\rm(iii)}] All of $A,\,B,\,C$ have three composition factors with the following structures,
\begin{eqnarray}
A:\ \ V_3\ \ \
\,\put(4,4){$\vector(3,1){24}$}\put(4,0){$\vector(3,-1){24}$}\put(31,10){$V_1$}
\put(31,-10){$V_2$}\hspace*{44pt}\ ,&\ \ \ \ \ \ \ \ B:\ \
\put(16,-8){$\vector(3,1){24}$}\put(16,10){$\vector(3,-1){24}$}\put(1,10){$V_1$}
\put(1,-10){$V_2$}\hspace*{42pt}V_3\ , \ \ \ \ \ \ \ \ &C:\ \ V_1\to
V_3\to V_2\ ,
\end{eqnarray}where $V_1=\F v_{0,-1},\,V_2=\F v_{0,-2},\,V_3={\rm
span\,}\{v_{k,l}\,|\,(k,l)\not=(0,-1),(0,-2)\}$,
\end{itemize}
\end{theo}
\vskip7pt

\cl{\bf\S2. \ Some technical lemmas
}\setcounter{section}{2}\setcounter{equation}{0} \vs{5pt}

It is well-known that the  Virasoro algebra $Vir$ is the Lie algebra
with basis $\{L_i,c\,|\,i\in\Z\}$ satisfying the relations
\begin{equation}\label{Vir}
[L_i,L_j]=(j-i)L_{i+j}+\frac{i^{3}-i}{12}\delta_{i,-j}c,\ \ \ \
[L_i,c]=0\mbox{\ \ \ for \ \ }i,j\in \Z.
\end{equation}
 A $Vir$-module of the intermediate series must be one of $A_{a,b},
A(a), B(a), a, b\in\C$, or one of the quotient submodules, where
$A_{a,b}, A(a), B(a)$ all have a basis $\{x_k\,|\,k\in\Z\}$ such
that $c$ acts trivially and
\begin{eqnarray}\label{condition}
\!\!\!\!\!\!\!\!\!\!\!\!&\!\!\!\!\!\!\!\!\!\!\!\!\!\!\!& A_{a,b}:
L_ix_k=(a+k+bi)x_{i+k},
\\\label{condition+}\!\!\!\!\!\!\!\!\!\!\!\!&\!\!\!\!\!\!\!\!\!\!\!\!\!\!\!&
A(a): L_ix_k=(i+k)x_{i+k}\ \ (k\neq0),
\ \ \ \ \ \ \ \ L_ix_0=i(i+a)x_i,
\\\label{condition++}\!\!\!\!\!\!\!\!\!\!\!\!&\!\!\!\!\!\!\!\!\!\!\!\!\!\!\!&
B(a): L_ix_k=kx_{i+k}\ \ (k\neq-i),
\ \ \ \ \ \ \ \ \ \ \ L_ix_{-i}=-i(i+a)x_0,
\end{eqnarray}
for $i,k\in\Z.$ We have
\begin{equation} A_{a,b} {\rm \ is\  simple}
\Longleftrightarrow a \notin\Z {\rm \ or\ } a \in\Z, b\neq0,1 {\rm \
and\ } A_{a,1}\cong A_{a,0} {\rm \ if\ } a\notin\Z.\label{isomo}
\end{equation}

Now consider the Lie algebra $\LL$ in (\ref{def1}). Denote
\begin{equation*}\mbox{
$\LL_0={\rm span}\{L_{i,0}\,|\,i\in\Z\}$ \ \ and \ $\LL_1={\rm
span}\{L_{i,i}\,|\,i\in\Z\}$. }\end{equation*} From (\ref{def1}), we
have
\begin{equation}\label{REL}\mbox{$[L_{i,0},L_{j,0}]=(j-i)L_{i+j,0}$ \ \ and
\ \ $[L_{i,i},L_{j,j}]=(j-i)L_{i+j,i+j}$ for all
$i,j\in\Z.$}\end{equation} This shows that both $\LL_0$ and $\LL_1$
are subalgebras of $\LL$ isomorphic to the centerless Virasoro
algebra (cf.~(\ref{Vir})). For any $j\in\Z,$ we denote
\begin{equation*}\mbox{$ V_j=\raisebox{-7pt}{$\stackrel{\dis \oplus}{\sc i\in\Z}$}\F v_{i,j}$.}\end{equation*}
Then $V_j$ is an $\LL_0$-module for all $j\in\Z$.

First we assume that all $V_j$ is an $\LL_0$-module of type
$A_{a,b}$ (later on, we shall determine all possible deformations,
cf.~(\ref{general})) . Thus we can assume
\begin{equation} L_{i,0}v_{k,j}=(a_j+k+b_ji)v_{i+k,j} \ {\rm \ for\
}\ i,k,j\in\Z \ {\rm \ and\ some}\ a_j,b_j\in\F.\label{conclu}
\end{equation}
Since $L_{0,0}v_{k,j}=(a_j+k)v_{k,j}$ and $L_{i',j'}v_{k,j}\subset
V_{j+j'}$ for any $i',j',k,j\in\Z$, applying $[L_{0,0},
L_{i',j'}]=i'L_{i',j'}$ to $v_{k,j}$, we can deduce that for any
$j\in\Z$, $a_j=a$ is a constant. Thus (\ref{conclu}) becomes
\begin{equation} L_{i,0}v_{k,j}=(a+k+b_ji)v_{i+k,j} \ {\rm \ for\
}\ i,k,j\in\Z.\label{conclu1}
\end{equation}
Analogously, we denote $\tilde V_1=\oplus_{i\in\Z}\F v_{i,i}$. Then
it is an $\LL_1$-module. Again, we first suppose $\tilde V_1$ is an
$\LL_1$-module of type $A_{a,c}$. By rescaling basis
$v_{i,i},\,i\in\Z$ (this rescaling does not affect (\ref{conclu1})
since for each given $j$, we can rescaling $v_{k,j}$ for all $k$
accordingly), we can suppose
\begin{equation} L_{i,i}v_{j,j}=(a+j+ci)v_{i+j,i+j} \ {\rm \ for\
}\ i,j\in\Z \ {\rm \ and\ some}\ c\in\F.\label{conclu2}
\end{equation}
Since $V=\oplus_{i,j\in \Z}\F v_{i,j}$ is an $\LL$-module of the
intermediate series, we can suppose
\begin{equation}\label{conclu3}L_{i,j}v_{k,l}=a^{i,j}_{k,l}v_{i+k,j+l} \ {\rm \ for\ }\
i,j,k,l\in\Z \ {\rm \ and\ some\ }\
a^{i,j}_{k,l}\in\F.\end{equation} Our aim is to prove that
\begin{equation}\label{To-Prove}\mbox{
$a^{i,j}_{k,l}=(j+1)(a+k)+(b_0-l)i$ \ \ for \
$i,j,k,l\in\Z$.}\end{equation} So we start with the case when $V_j$
has a basis satisfying  (\ref{conclu1}) for all $j$ and $\tilde V_1$
satisfying (\ref{conclu3}).
\begin{lemm}\label{lemm} \rm $b_s=b_0-s$ for all $s\in\Z.$
\end{lemm}
\noindent{\it Proof.~} Fix $s\in\Z$. For $i,j,k\in\Z,$ using
formulas
\begin{eqnarray}\label{FOr1}&&[L_{i,0},[L_{j,0},L_{k,s}]]=(j+k-(s+1)i)(k-(s+1)j)L_{i+j+k,s},\\
&&\label{FOr2}
[L_{i+j,0},L_{k,s}]=(k-(i+j)(s+1))L_{i+j+k,s},\end{eqnarray} we
obtain
\begin{eqnarray}\label{equ}
\!\!\!\!\!\!\!\!\!\!\!\!&\!\!\!\!\!\!\!\!\!\!\!\!\!\!\!&
(k-(i+j)(s+1))(L_{i,0}L_{j,0}L_{k,s}-L_{i,0}L_{k,s}L_{j,0}-L_{j,0}L_{k,s}L_{i,0}+L_{k,s}L_{j,0}L_{i,0})
\nonumber\\\!\!\!\!\!\!\!\!\!\!\!\!&\!\!\!\!\!\!\!\!\!\!\!\!\!\!\!&
\,\,\,\,\,\,\,\,\,\,\,\,\,\,\,\,\,\,\,\,\,
-(j+k-(s+1)i)(k-(s+1)j)(L_{i+j,0}L_{k,s}-L_{k,s}L_{i+j,0})=0.
\end{eqnarray}
Now applying (\ref{equ}) to
\begin{equation}\label{v-0}v_{p,0},\end{equation} by (\ref{conclu1}) and
(\ref{conclu3}), and by comparing the coefficients of
$v_{i+j+k+p,s}$, we obtain
\begin{eqnarray}\label{four}
\!\!\!\!\!\!\!\!\!\!\!\!&\!\!\!\!\!\!\!\!\!\!\!\!\!\!\!&((i\!+\!j)(s\!+\!1)\!-\!k)
(a\!+\!p\!+\!jb_0)(a\!+\!j\!+\!k\!+\!p\!+\!ib_s)a^{k,s}_{j+p,0}
\nonumber\\\!\!\!\!\!\!\!\!\!\!\!\!&\!\!\!\!\!\!\!\!\!\!\!\!\!\!\!&
\ \ \ \ \ \ \!+
((k\!-\!(i\!+\!j)(s\!+\!1))(a\!+\!p\!+\!k\!+\!jb_s)(a\!+\!p\!+\!k\!+\!j\!+\!ib_s)
\nonumber\\\!\!\!\!\!\!\!\!\!\!\!\!&\!\!\!\!\!\!\!\!\!\!\!\!\!\!\!&
\ \ \ \ \ \
\!-(j\!+\!k\!-\!(s\!+\!1)i)(k\!-\!j(s\!+\!1))(a\!+\!p\!+\!k\!+\!b_s(i\!+\!j))a^{k,s}_{p,0}
\nonumber\\\!\!\!\!\!\!\!\!\!\!\!\!&\!\!\!\!\!\!\!\!\!\!\!\!\!\!\!&
\ \ \ \ \ \ \!+
((k\!-\!(i\!+\!j)(s\!+\!1))(a\!+\!p\!+\!ib_0)(a\!+\!p\!+\!i\!+\!jb_0)
\nonumber\\\!\!\!\!\!\!\!\!\!\!\!\!&\!\!\!\!\!\!\!\!\!\!\!\!\!\!\!&
\ \ \ \ \ \
\!+(j\!+\!k\!-\!(s\!+\!1)i)(k\!-\!j(s\!+\!1))(a\!+\!p\!+\!b_0(i\!+\!j))a^{k,s}_{p+i+j,0}
\nonumber\\\!\!\!\!\!\!\!\!\!\!\!\!&\!\!\!\!\!\!\!\!\!\!\!\!\!\!\!&
\ \ \ \ \ \
\!+((i\!+\!j)(s\!+\!1)\!-\!k)(a\!+\!p\!+\!ib_0)(a\!+\!i\!+\!k\!+\!p\!+\!jb_s)a^{k,s}_{i+p,0}=0.
\end{eqnarray}
 Now in
(\ref{four}), by doing the following:
\begin{itemize}\parskip-3pt\item[$(1)$] replacing $j$ by $i$, and $p$ by
$p-i$; \item[$(2)$] replacing $j$ by $-i$; \item[$(3)$] replacing
both $i$ and $j$ by $-i$, and $p$ by $p+i$, \end{itemize} we obtain
the following three equations concerning $a^{k,s}_{p-i,0}$,
$a^{k,s}_{p,0}$, $a^{k,s}_{p+i,0}$ correspondingly:
\begin{eqnarray}\label{abc1}
\!\!\!\!\!\!\!\!\!\!\!\!&\!\!\!\!\!\!\!\!\!\!\!\!\!\!\!&((k\!-\!2i(s\!+\!1))(a\!+\!p\!+\!k\!+\!i(b_s\!-\!1))(a\!+\!k\!+\!p\!+\!ib_s)
\nonumber\\\!\!\!\!\!\!\!\!\!\!\!\!&\!\!\!\!\!\!\!\!\!\!\!\!\!\!\!&
\ \ \ \ \ \ \ \ \ \ -
(k\!-\!is)(k\!-\!(s\!+\!1)i)(a\!+\!p\!+\!k\!+\!i(2b_s\!-\!1)))a^{k,s}_{p-i,0}
\nonumber\\\!\!\!\!\!\!\!\!\!\!\!\!&\!\!\!\!\!\!\!\!\!\!\!\!\!\!\!&
\ \ \ \ \ \ \ \ \ \
+2(2i(s\!+\!1)\!-\!k)(a\!+\!p\!+\!(b_0\!-\!1)i)(a\!+\!p\!+\!k\!+\!b_si)a^{k,s}_{p,0}
\nonumber\\\!\!\!\!\!\!\!\!\!\!\!\!&\!\!\!\!\!\!\!\!\!\!\!\!\!\!\!&
\ \ \ \ \ \ \ \ \ \ +
((k\!-\!2i(s\!+\!1))(a\!+\!p\!+\!ib_0)(a\!+\!p\!+\!i(b_0\!-\!1))
\nonumber\\\!\!\!\!\!\!\!\!\!\!\!\!&\!\!\!\!\!\!\!\!\!\!\!\!\!\!\!&
\ \ \ \ \ \ \ \ \ \
+(k\!-\!(s\!+\!1)i)(k\!-\!is)(a\!+\!p\!+\!(2b_0\!-\!1)i))a^{k,s}_{p+i,0}
=0,
\end{eqnarray}
\begin{eqnarray}
\label{abc2}
\!\!\!\!\!\!\!\!\!\!\!\!&\!\!\!\!\!\!\!\!\!\!\!\!\!\!\!&-k(a\!+\!p\!-\!b_0i)(a\!+\!p\!+\!k\!+\!i(b_s\!-\!1))a^{k,s}_{p-i,0}\!+\!
k((a\!+\!p\!+\!k\!-\!b_si)\times
\nonumber\\\!\!\!\!\!\!\!\!\!\!\!\!&\!\!\!\!\!\!\!\!\!\!\!\!\!\!\!&
\ \ \ \ \ \ \ \ \ \ \
\times(a\!+\!p\!+\!k\!+\!i(b_s\!-\!1))\!+\!(a\!+\!p\!+\!b_0i)(a\!+\!p\!-\!(b_0\!-\!1)i)
\nonumber\\\!\!\!\!\!\!\!\!\!\!\!\!&\!\!\!\!\!\!\!\!\!\!\!\!\!\!\!&
\ \ \ \ \ \ \ \ \ \ \
-(k\!-\!i(s\!+\!2))(k\!+\!i(s\!+\!1)))a^{k,s}_{p,0}\!-k\!(a\!+\!p\!+\!b_0i)\times
\nonumber\\\!\!\!\!\!\!\!\!\!\!\!\!&\!\!\!\!\!\!\!\!\!\!\!\!\!\!\!&
\ \ \ \ \ \ \ \ \ \ \
\times(a\!+\!p\!+\!k\!-\!i(b_s\!-\!1))a^{k,s}_{p+i,0} =0,
\\[4pt] 
\label{abc3}
\!\!\!\!\!\!\!\!\!\!\!\!&\!\!\!\!\!\!\!\!\!\!\!\!\!\!\!&\
((k\!+\!2i(s\!+\!1))(a\!+\!p\!-\!(b_0\!-\!1)i)(a\!+\!p\!-\!ib_0)
\nonumber\\\!\!\!\!\!\!\!\!\!\!\!\!&\!\!\!\!\!\!\!\!\!\!\!\!\!\!\!&
\ \ \ \ \ \ \ \ \ \ \
+(k\!+\!is)(k\!+\!(s\!+\!1)i)(a\!+\!p\!-\!i(2b_0\!-\!1)))a^{k,s}_{p-i,0}
\nonumber\\\!\!\!\!\!\!\!\!\!\!\!\!&\!\!\!\!\!\!\!\!\!\!\!\!\!\!\!&
\ \ \ \ \ \ \ \ \ \ \
-2(2i(s\!+\!1)\!+\!k)(a\!+\!p\!-\!(b_0\!-\!1)i)(a\!+\!p\!+\!k\!-\!b_si)a^{k,s}_{p,0}
\nonumber\\\!\!\!\!\!\!\!\!\!\!\!\!&\!\!\!\!\!\!\!\!\!\!\!\!\!\!\!&
\ \ \ \ \ \ \ \ \ \ \ +
((k\!+\!2i(s\!+\!1))(a\!+\!p\!+\!k\!-\!i(b_s\!-\!1))(a\!+\!p\!+\!k\!-\!ib_s)
\nonumber\\\!\!\!\!\!\!\!\!\!\!\!\!&\!\!\!\!\!\!\!\!\!\!\!\!\!\!\!&
\ \ \ \ \ \ \ \ \ \ \
-(k\!+\!(s\!+\!1)i)(k\!+\!is)(a\!+\!p\!+\!k\!-\!(2b_s\!-\!1)i)a^{k,s}_{p+i,0}
=0.
\end{eqnarray}
Denote the matrix of the coefficients of $a^{k,s}_{p-i,0}$,
$a^{k,s}_{p,0}$, $a^{k,s}_{p+i,0}$ in (\ref{abc1}), (\ref{abc2}) and
(\ref{abc3}) by $A=(c_{ij})_{3\times3}$ and let $\Delta=\det\, A$.
Observe that if we denote $c_{31}$ by $f_1(i)$, $c_{32}$ by
$f_2(i)$, $c_{33}$ by $g_2(i)$ and $c_{23}$ by $g_1(i)$ $($all
regarding as polynomials on $i),$ then (\ref{abc1})--(\ref{abc3})
become
\begin{eqnarray}\label{equ1}
&&
g_2(-i)a^{k,s}_{p-i,0}+f_2(-i)a^{k,s}_{p,0}+f_1(-i)a^{k,s}_{p+i,0}=0,\\
\label{equ2}
&&g_1(-i)a^{k,s}_{p-i,0}+c_{22}a^{k,s}_{p,0}+g_1(i)a^{k,s}_{p+i,0}=0,\\
\label{equ3}
&&f_1(i)a^{k,s}_{p-i,0}+f_2(i)a^{k,s}_{p,0}+g_2(i)a^{k,s}_{p+i,0}=0.
\end{eqnarray}
Note from definition (\ref{conclu3}) that for a fixed $s$, one can
easily prove $a^{k,s}_{p,0}\neq0$ for almost all $k,p\in\Z,$ we have
$\Delta=0$ (for almost all $i,k,p$, thus for all $i,k,p$ since
$\Delta$ is a polynomial on $i,k,p$). By a little lengthy
calculation $($or simply using Mathematica to compute$)$, we obtain
\begin{eqnarray}\label{det}
\!\!\!\!\!\!\!\!\!\!\!\!&\!\!\!\!\!\!\!\!\!\!\!\!\!\!\!&
\Delta=i^6k(s\!-\!b_0\!+\!b_s)(1\!+\!s\!-\!b_0\!+\!b_s)(\Delta_{0}\!+\!\Delta_{1}k\!+\!\Delta_{2}k^2),
\\\!\!\!\!\!\!\!\!\!\!\!\!&\!\!\!\!\!\!\!\!\!\!\!\!\!\!\!&
\Delta_{0}=4i^2(1\!+\!s)^2(-1\!+\!b_0\!+\!b_s)(2b_0\!+\!3sb_0\!+\!s^2b_0\!+\!b_0^2\!-\!b_0^3\!+\!2b_s\!-\!sb_s
\nonumber\\\!\!\!\!\!\!\!\!\!\!\!\!&\!\!\!\!\!\!\!\!\!\!\!\!\!\!\!&
\ \ \ \ \ \ \ \
-s^2b_s\!-\!2b_0b_s\!-\!2sb_0b_s\!-\!b_0^2b_s\!-\!3b_s^2\!+\!b_0b_s^2\!+\!b_s^3),
\nonumber\\\!\!\!\!\!\!\!\!\!\!\!\!&\!\!\!\!\!\!\!\!\!\!\!\!\!\!\!&
\Delta_{1}=2(a\!+\!p)(1\!+\!s)(2\!+\!s)(-2\!+\!s\!+\!s^2\!+\!7b_0\!+\!2sb_0\!-\!3b_0^2\!+\!5b_s\!-\!2sb_s
\!-\!6b_0b_s\!-\!3b_s^2),
\nonumber\\\!\!\!\!\!\!\!\!\!\!\!\!&\!\!\!\!\!\!\!\!\!\!\!\!\!\!\!&
\Delta_{2}=-4\!+\!24b_0\!-\!19b_0^2\!+\!2b_0^3\!+\!b_0^4\!+\!16b_s\!-\!34b_0b_s\!+\!8b_0^2b_s\!+\!2b_0^3b_s\!-\!19b_s^2
\nonumber\\\!\!\!\!\!\!\!\!\!\!\!\!&\!\!\!\!\!\!\!\!\!\!\!\!\!\!\!&
\ \ \ \ \ \ \ \
+14b_0b_s^2\!+\!8b_s^3\!-\!2b_0b_s^3\!-\!b_s^4\!+\!2(-2\!+\!18b_0\!-\!11b_0^2\!+\!b_0^3\!+\!4b_s\!-\!15b_0b_s
\nonumber\\\!\!\!\!\!\!\!\!\!\!\!\!&\!\!\!\!\!\!\!\!\!\!\!\!\!\!\!&
\ \ \ \ \ \ \ \
+3b_0^2b_s\!-\!4b_s^2\!+\!3b_0b_s^2\!+\!b_s^3)s\!+\!(3\!+\!16b_0\!-\!6b_0^2\!-\!4b_s\!-\!6b_0b_s)s^2
\!+\!2(2\!+\!b_0\!-\!b_s)s^3\!+\!s^4.\nonumber
\end{eqnarray} Thus we have
\begin{eqnarray}\label{b-con1}
&&(s-b_0+b_s)(1+s-b_0+b_s)=0, \ \ \rm or\\
\label{b-con2}&& \Delta_0=\D_1=\D_2=0.\end{eqnarray} Note that if we
replace $v_{p,0}$ (in (\ref{v-0})) and $s$ by $v_{p,s}$ and $-s$
respectively (and accordingly, replace $a^{k,s}_{p,0}$ by
$a^{k,-s}_{p,s}$), this amounts to change the data $(b_0,b_s,s)$ to
$(b_s,b_0,-s)$ in (\ref{det}). Thus we have
\begin{eqnarray}\label{b-con1+}
&&(-s-b_s+b_0)(1-s-b_s+b_0)=0, \ \ \rm or\\
\label{b-con2+}&& \Delta'_0=\D'_1=\D'_2=0,\end{eqnarray} where
$\D'_i$ is obtained from $\D_i$ by changing the data $(b_0,b_s,s)$
to $(b_s,b_0,-s)$. We want to prove $b_s=b_0-s$. Suppose this is not
the case. Then we have only three possible cases from
(\ref{b-con1})--(\ref{b-con2+}):
\begin{itemize}\parskip-3pt\item[(i)] $b_s=b_0-1-s$ and
$\D'_0=\D'_1=\D'_2=0$;
\item[(ii)] $\D_0=\D_1=\D_2=0$ and $b_s=b_0+1-s$;
 \item[(iii)] $\D_i=\D'_i=0,\,i=0,1,2$.\end{itemize}
It is a little lengthy but straightforward to verify that these
three cases are impossible. This completes the proof of the
lemma.\hfill$\Box$\vskip5pt

Using (\ref{conclu1}), (\ref{conclu2}) and Lemma \ref{lemm}, we
obtain
\begin{equation}\label{coeff}a^{i,i}_{j,j}=a+j+ci,\ \ \ a^{i,0}_{k,l}=a+k+(b_0-l)i\ \ \
 {\rm \ for\ \ }\ i,j,k,l\in\Z.
\end{equation}
\begin{lemm}\label{lemm2} \rm $((b_0-q)k+(s+1)(a+p))a^{k,s}_{p+i,q}=((b_0-q)k+(s+1)(a+p+i))a^{k,s}_{p,q}$ for
 $i,k,s,p,q\in\Z$.
\end{lemm}
\noindent{\it Proof.~} Applying (\ref{equ}) to $v_{p,q}$, by
(\ref{conclu1}) and (\ref{conclu3}), comparing the coefficients of
$v_{i+j+k+p,s+q}$, we obtain \begin{eqnarray}\label{non}
\!\!\!\!\!\!\!\!\!\!\!\!&\!\!\!\!\!\!\!\!\!\!\!\!\!\!\!&
((i\!+\!j)(s\!+\!1)\!-\!k)(a\!+\!p\!+\!j(b_0\!-\!q))(a\!+\!j\!+\!k\!+\!p\!+\!i(b_0\!-\!q\!-\!s))a^{k,s}_{j+p,q}
\nonumber\\\!\!\!\!\!\!\!\!\!\!\!\!&\!\!\!\!\!\!\!\!\!\!\!\!\!\!\!&
\ \ \ \ +
((k\!-\!(i\!+\!j)(s\!+\!1))(a\!+\!p\!+\!k\!+\!j(b_0\!-\!q\!-\!s))(a\!+\!p\!+\!k\!+\!j\!+\!i(b_0\!-\!q\!-\!s))
\nonumber\\\!\!\!\!\!\!\!\!\!\!\!\!&\!\!\!\!\!\!\!\!\!\!\!\!\!\!\!&
\ \ \ \
-(j\!+\!k\!-\!(s\!+\!1)i)(k\!-\!j(s\!+\!1))(a\!+\!p\!+\!k\!+\!(b_0\!-\!q\!-\!s)(i\!+\!j))a^{k,s}_{p,q}
\nonumber\\\!\!\!\!\!\!\!\!\!\!\!\!&\!\!\!\!\!\!\!\!\!\!\!\!\!\!\!&
\ \ \ \ +
((k\!-\!(i\!+\!j)(s\!+\!1))(a\!+\!p\!+\!i(b_0\!-\!q))(a\!+\!p\!+\!i\!+\!j(b_0\!-\!q))
\nonumber\\\!\!\!\!\!\!\!\!\!\!\!\!&\!\!\!\!\!\!\!\!\!\!\!\!\!\!\!&
\ \ \ \
+(j\!+\!k\!-\!(s\!+\!1)i)(k\!-\!j(s\!+\!1))(a\!+\!p\!+\!(b_0\!-\!q)(i\!+\!j))a^{k,s}_{p+i+j,q}
\nonumber\\\!\!\!\!\!\!\!\!\!\!\!\!&\!\!\!\!\!\!\!\!\!\!\!\!\!\!\!&
\ \ \ \
+((i\!+\!j)(s\!+\!1)\!-\!k)(a\!+\!p\!+\!i(b_0\!-\!q))(a\!+\!i\!+\!k\!+\!p\!+\!j(b_0\!-\!q\!-\!s))a^{k,s}_{i+p,q}=0.
\end{eqnarray}
The same arguments after (\ref{four}) can give three equations
concerning $a^{k,s}_{p-i,q}$, $a^{k,s}_{p,q}$, $a^{k,s}_{p+i,q}$.
Solving them gives
\begin{eqnarray}\label{add}\!\!\!\!\!\!\!\!\!\!\!\!&\!\!\!\!\!\!\!\!\!\!\!\!\!\!\!&
((b_0-q)k+(s+1)(a+p))a^{k,s}_{p+i,q}=((b_0-q)k+(s+1)(a+p+i))a^{k,s}_{p,q},
\end{eqnarray}
for  $k,s,i,p,q\in\Z$. Now the lemma follows.\hfill$\Box$\vskip5pt
\begin{lemm}\label{lemm3} \rm $c=a+b_0$.
\end{lemm}
\noindent{\it Proof.~} For $i,j,k,l\in\Z,$ applying
$[L_{i,i},L_{j,j}]=(j-i)L_{i+j,i+j}$
to $v_{k,l}$, by (\ref{conclu2}) and (\ref{conclu3}), comparing the
coefficients of $v_{i+j+k,i+j+l}$, we obtain
\begin{equation}\label{sym}
a^{j,j}_{k,l}a^{i,i}_{j+k,j+l}-a^{i,i}_{k,l}a^{j,j}_{i+k,i+l}-(j-i)a^{i+j,i+j}_{k,l}=0.
\end{equation}
Changing data $(s,p,i)$ in (\ref{add}) to $(k,q,i-q)$ gives
\begin{equation}\label{ss}(q+b_0k+(k+1)a)a^{k,k}_{i,q}=((b_0-q)k+(k+1)(a+i))(a+q+ck).\end{equation}
Then (\ref{sym}) is equivalent to the following equation
\begin{eqnarray}\label{hard}
\!\!\!\!\!\!\!\!\!\!\!\!&\!\!\!\!\!\!\!\!\!\!\!\!\!\!\!&(a\!+\!b_0\!-\!c)i
j(j\!-\!i)(k\!-\!l)((a\!+\!b_0)c(i\!+\!j)\!+\!(a\!+\!l)(a\!+\!b_0\!+\!c\!-\!1))
A=0,
\\\!\!\!\!\!\!\!\!\!\!\!\!&\!\!\!\!\!\!\!\!\!\!\!\!\!\!\!&
 \ A=(a\!+\!b_0)(l\!+\!j)(a\!+\!b_0\!+\!k\!-\!l)i^2\!+\!(2a^2\!+\!b_0j\!+\!b_0^2j\!+\!b_0^2j^2\!+\!b_0k\!+\!2b_0jk
\nonumber\\\!\!\!\!\!\!\!\!\!\!\!\!&\!\!\!\!\!\!\!\!\!\!\!\!\!\!\!&\
\ \ \ \ \ \ +b_0j^2k\!+\!b_0l\!-\!j
l\!-\!b_0j^2l\!+\!kl\!+\!jkl\!-\!l^2\!-\!jl^2\!+\!2ab_0\!+\!3a^2j\!+\!4abj
\nonumber\\\!\!\!\!\!\!\!\!\!\!\!\!&\!\!\!\!\!\!\!\!\!\!\!\!\!\!\!&\
\ \ \ \ \ \ +a^2j^2\!+\!2abj^2\!+\!2ak\!+\!3ajk\!+\!aj^2k\!-\!a j
l\!-\!aj^2l)i
\nonumber\\\!\!\!\!\!\!\!\!\!\!\!\!&\!\!\!\!\!\!\!\!\!\!\!\!\!\!\!&\
\ \ \ \ \ \
+(a\!+\!aj\!+\!b_0j\!+\!l)(a\!+\!aj\!+\!b_0j\!+\!k\!+\!jk\!-\!j
l).\nonumber
\end{eqnarray}
Since (\ref{hard}) is a polynomial on $i,j,k,l$, we must have
$c=a+b_0$ or $(a+b_0)c=a+b_0+c-1=0$. Thus $c=a+b_0$, or
$a+b_0=0,c=1$, or $a+b_0=1,c=0$. If the first case occurs, then we
have the lemma. Assume that one of the last two cases occurs. In
case $a\notin\Z$, by (\ref{isomo}), we can choose $c$ to be 0 or 1,
so that $c=a+b_0$. If $a\in\Z$, we can rescale basis
$\{v_{ii}\,|\,i\in\Z\}$ in (\ref{conclu2}) so that $c$ can be chosen
either one of $0$ or $1$ (in this case, (\ref{conclu2}) holds only
when $j,i+j\ne0$, but we can first assume it holds for all $i,j$,
then consider all possible deformations (which will be done in the
last part of the next section)). Thus without loss of generality, we
can always suppose $c=a+b_0$.\hfill$\Box$\vskip5pt

Therefore by equation (\ref{conclu2}), we obtain
\begin{equation}\label{double}
a^{i,i}_{k,k}=a+k+(a+b_0)i\ \  {\rm \ for\ all}\ i,k\in\Z.
\end{equation}
Replacing $p+i$ by $i$ and then letting $p=q$ in (\ref{add}), we
obtain
\begin{equation}\label{near}((b_0-q)k+(s+1)(a+q))a^{k,s}_{i,q}=((s+1)(a+i)+(b_0-q)k)
a^{k,s}_{q,q}\ {\rm \ for\ any}\ k,s,i,q\in\Z. \end{equation}

\begin{lemm}\label{lemm4}
\rm For any $p,q,k,s\in\Z$, we have \begin{eqnarray}\label{replace}
\!\!\!\!\!\!\!\!\!\!\!\!&\!\!\!\!\!\!\!\!\!\!\!\!\!\!\!&
((b_0-q)k+(s+1)(a+q))(k-(s+1)(p-k))a^{p,s}_{q,q}
\nonumber\\\!\!\!\!\!\!\!\!\!\!\!\!&\!\!\!\!\!\!\!\!\!\!\!\!\!\!\!&
=(((b_0-q)k+(s+1)(a+q))(k-s(p-k))
\nonumber\\\!\!\!\!\!\!\!\!\!\!\!\!&\!\!\!\!\!\!\!\!\!\!\!\!\!\!\!&
\ \ \ -(a+q+(b_0-q)(p-k))(s+1)(p-k))a^{k,s}_{q,q}.\end{eqnarray}
\end{lemm}
\noindent{\it Proof.~} For any $i,k,s,q\in\Z$, applying
$[L_{i,0},L_{k,s}]=(k-(s+1)i)L_{i+k,s}$ to $v_{q,q}$, 
by (\ref{coeff}) and (\ref{double}), comparing the coefficients of
$v_{i+k+q,s+q}$, we have
\begin{eqnarray}\label{b}
\!\!\!\!\!\!\!\!\!\!\!\!&\!\!\!\!\!\!\!\!\!\!\!\!\!\!\!&
((b_0-q)k+(s+1)(a+q))(k-s i)-(a+q+(b_0-q)i)(s+1)i)a^{k,s}_{q,q}
\nonumber\\\!\!\!\!\!\!\!\!\!\!\!\!&\!\!\!\!\!\!\!\!\!\!\!\!\!\!\!&
=((b_0-q)k+(s+1)(a+q))(k-(s+1)i)a^{i+k,s}_{q,q}.
\end{eqnarray}
Replacing $i+k$ by $p$ in (\ref{b}), we immediately obtain that the
lemma holds.\hfill$\Box$\vskip8pt

\cl{\bf\S3. \ Proof of Theorem \ref{theo}
}\setcounter{section}{3}\setcounter{equation}{0} \vs{5pt} Now we can
prove the main results of this paper.

\ni{\it Proof of Theorem \ref{theo}.} For any $k,s,i,t\in\Z$, by
(\ref{coeff}), (\ref{double}), (\ref{near}), and Lemma \ref{lemm4},
we obtain
\begin{eqnarray}\label{equation}
\!\!\!\!\!\!\!\!\!\!\!\!&\!\!\!\!\!\!\!\!\!\!\!\!\!\!\!&
((b_0\!-\!t)s\!+\!(s\!+\!1)(a\!+\!t))(s\!-\!(s\!+\!1)(k\!-\!s))((b_0\!-\!t)k\!+\!(s\!+\!1)(a\!+\!t))a^{k,s}_{i,t}
\nonumber\\\!\!\!\!\!\!\!\!\!\!\!\!&\!\!\!\!\!\!\!\!\!\!\!\!\!\!\!&
=((s\!+\!1)(a\!+\!i)\!+\!(b_0\!-\!t)k)((b_0\!-\!t)s\!+\!(s\!+\!1)(a\!+\!t))(s\!-\!(s\!+\!1)(k\!-\!s))
a^{k,s}_{t,t}
\nonumber\\\!\!\!\!\!\!\!\!\!\!\!\!&\!\!\!\!\!\!\!\!\!\!\!\!\!\!\!&
=((s\!+\!1)(a\!+\!i)\!+\!(b_0\!-\!t)k)(((b_0\!-\!t)s\!+\!(s\!+\!1)(a\!+\!t))(s\!-\!s(k\!-\!s))
\nonumber\\\!\!\!\!\!\!\!\!\!\!\!\!&\!\!\!\!\!\!\!\!\!\!\!\!\!\!\!&
\ \ \
-(a\!+\!t\!+\!(b_0\!-\!t)(k\!-\!s))(s\!+\!1)(k\!-\!s))a^{s,s}_{t,t}
\nonumber\\\!\!\!\!\!\!\!\!\!\!\!\!&\!\!\!\!\!\!\!\!\!\!\!\!\!\!\!&
=((s\!+\!1)(a\!+\!i)\!+\!(b_0\!-\!t)k)(s\!-\!(s\!+\!1)(k\!-\!s))\times
\nonumber\\\!\!\!\!\!\!\!\!\!\!\!\!&\!\!\!\!\!\!\!\!\!\!\!\!\!\!\!&
\ \ \
\times((b_0\!-\!t)k\!+\!(s\!+\!1)(a\!+\!t))((b_0\!-\!t)s\!+\!(s\!+\!1)(a\!+\!t)).
\end{eqnarray}
Therefore, we obtain
\begin{eqnarray}\label{fin}
\!\!\!\!\!\!\!\!\!\!\!\!&\!\!\!\!\!\!\!\!\!\!\!\!\!\!\!&
((b_0\!-\!t)k\!+\!(s\!+\!1)(a\!+\!t))((b_0\!-\!t)s\!+\!(s\!+\!1)(a\!+\!t))(s\!-\!(s\!+\!1)(k\!-\!s))\times
\nonumber\\\!\!\!\!\!\!\!\!\!\!\!\!&\!\!\!\!\!\!\!\!\!\!\!\!\!\!\!&
\times(a^{k,s}_{i,t}\!-\!((s\!+\!1)(a\!+\!i)\!+\!(b_0\!-\!t)k))=0.
\end{eqnarray}
From (\ref{fin}), we obtain
\begin{equation}\label{a-i-j-k-l}
a_{i,t}^{k,s}=(s+1)(a+i)+(b_0-t)k,
\end{equation}
if the $4$-tuple $(k,s,i,t)$ satisfies the condition
\begin{equation}\label{condition-a-i-j-k-l}
((b_0-t)k+(s+1)(a+t))((b_0-t)s+(s+1)(a+t))(s-(s+1)(k-s))\ne0.
\end{equation}
We want to prove (\ref{a-i-j-k-l}) holds for all $4$-tuple
$(k,s,i,t)$. By (\ref{coeff}), we have

\begin{equation}\label{a-i-j-k-0}
a^{k,0}_{i,t}=a+i+(b_0-t)k \mbox{ \ \ for any\ }i,k,t\in\Z.
\end{equation}
So we can suppose $s\ne0$. Applying
$[L_{k',s'},L_{k-k',s-s'}]=((s'+1)(k-k')-(s-s'+1)k')L_{k,s}$ to
$v_{i,t}$, and comparing the coefficients of $v_{k+i,s+t}$, we have
\begin{equation}\label{relation-a-i-j-k-l}
((s'+1)(k-k')-(s-s'+1)k')a_{i,t}^{k,s}=a_{k+i-k',t+s-s'}^{k',s'}a_{i,t}^{k-k',s-s'}
-a_{i+k',t+s'}^{k-k',s-s'}a_{i,t}^{k',s'}.
\end{equation}
By calculation, we know that if the following condition
\begin{equation}\label{condition-a-i-j}
(k,s)\not=(0,-2)  \mbox{ \ and\ } t\not=b_0
\end{equation}
holds then we can choose $(k',s')$ such that the $4$-tuples
$$(k',s',k+i-k',t+s-s'),\,(k-k',s-s',i,t),\,
(k-k',s-s',i+k',t+s'),\,(k',s',i,t)$$ satisfy condition
(\ref{condition-a-i-j-k-l}) and $(s'+1)(k-k')-(s-s'+1)k'\ne0$. Thus
under condition (\ref{condition-a-i-j}), (\ref{a-i-j-k-l}) follows
from (\ref{relation-a-i-j-k-l}). That is to say, we have
\begin{equation}\label{produce}
a_{i,t}^{k,s}=(s+1)(a+i)+(b_0-t)k {\rm \ \ for\ }k,s,i,t\in\Z,\
(k,s)\not=(0,-2) {\rm \ and\ } t\not=b_0.
\end{equation}
Now we only need to consider the following two cases.

{\it Case 1: $(k,s)=(0,-2)$ and $t\not=b_0$.}

For any $p,q,m,s\in\Z,$ applying
$[L_{p,q},L_{m,s}]=((q+1)m-(s+1)p)L_{p+m,q+s}$
 to $v_{i,t}$, comparing the coefficients of $v_{p+m+i,q+s+t}$,
we have
\begin{equation}\label{chief}
a^{p,q}_{i+m,s+t}a^{m,s}_{i,t}-a^{p,q}_{i,t}a^{m,s}_{i+p,t+q}=((q+1)m-(s+1)p)a^{p+m,q+s}_{i,t}.
\end{equation}
Set $(p,q,s)=(0,-2,0)$, $m\not=0$ and $t\not=b_0$ in (\ref{chief}),
by (\ref{coeff}) and (\ref{produce}), (\ref{chief}) turns into the
following equation
\begin{eqnarray}\label{sub1}
\!\!\!\!\!\!\!\!\!\!\!\!&\!\!\!\!\!\!\!\!\!\!\!\!\!\!\!&
(a\!+\!i\!+\!(b_0\!-\!t)m)a^{0,-2}_{i+m,t}\!-\!(a\!+\!i\!+\!(b_0\!-\!t\!+\!2)m)a^{0,-2}_{i,t}
\!=\!-m(-(a\!+\!i)\!+\!(b_0\!-\!t)m)).
\end{eqnarray}
By Lemma \ref{lemm2}, we can easily get
\begin{equation}\label{displace}
(a+m+i)a^{0,-2}_{i,t}=(a+i)a^{0,-2}_{m+i,t}.
\end{equation}
Therefore, multiplying (\ref{sub1}) by $a+i$, using
(\ref{displace}), we obtain
\begin{equation}\label{bak}
m(-(a+i)+(b_0-t)m)(a^{0,-2}_{i,t}+a+i)=0.
\end{equation}
Since by our assumption, we can choose suitable $m$ such that
$m(-(a+i)+(b_0-t)m)\not=0$, thus $a^{0,-2}_{i,t}=-(a+i)$ for
$i,t\in\Z$ and $t\not=b_0$. Therefore in this case we obtain
\begin{equation}\label{outcome}
a^{k,s}_{i,t}=(s+1)(a+i)+(b_0-t)k {\rm \ \ for\ }i,t\in\Z  {\rm \
and\ }(k,s)=(0,-2),\  t\not=b_0.
\end{equation}
Therefore, (\ref{outcome}) together with (\ref{produce}) shows that
\begin{equation}\label{outcome-t}
a^{k,s}_{i,t}=(s+1)(a+i)+(b_0-t)k {\rm \ \ for\ }k,s,i,t\in\Z  {\rm
\ and\ }t\not=b_0.
\end{equation}

{\it Case 2: $t=b_0$.}

In (\ref{chief}), we set $i+m=k$ and $t=b_0-s$, then we have
\begin{equation}\label{aks}
a^{p,q}_{k,b_0}a^{k-i,s}_{i,b_0-s}-a^{p,q}_{i,b_0-s}a^{k-i,s}_{i+p,b_0-s+q}=((q+1)(k-i)-(s+1)p)a^{p+k-i,q+s}_{i,b_0-s}.
\end{equation}
We can choose $s\not=0$ and $s\not=q$ such that $b_0-s\not=b_0$ and
$b_0-s+q\not=b_0$, by (\ref{outcome-t}) and by calculation, we have
the following equation
\begin{equation}\label{new1}
(a+i+(a+k)s)(a^{p,q}_{k,b_0}-(q+1)(a+k))=0.
\end{equation}
Since we can choose suitable $i$ and $s$ such that
$a+i+(a+k)s\not=0$, then we obtain
\begin{equation}\label{new2}
a^{p,q}_{k,b_0}=(q+1)(a+k)  {\rm \ \ for\ any} \ p,q,k\in\Z.
\end{equation}
This together with (\ref{outcome-t}) shows that
\begin{equation}\label{new3}
a^{k,s}_{i,t}=(s+1)(a+i)+(b_0-t)k  {\rm \ \ for\ any} \
k,s,i,t\in\Z.
\end{equation}
In a word, for any $i,j,k,l\in \Z,$ we have
$L_{i,j}v_{k,l}=((j+1)(a+k)+(b_0-l)i)v_{i+k,j+l}.$ Denote $b_0$ by
$b-1$, then we obtain that (\ref{A-a-b}) in Theorem \ref{theo}
holds.

Now in the following, we discuss the irreducibility of $A_{a,b}$.
Suppose $V'$ is a submodule of $A_{a,b}$ and
$0\not=v=\Sigma_{k=0}^{m}a_k v_{0,k}\in V'$, $a_k\not=0$ for all
$0\leqslant k\leqslant m$. By (\ref{A-a-b}), we have
\begin{equation}\label{action}
\rho\cdot v:=L_{i,0}L_{-i,0}\cdot
v=\mbox{$\sum\limits_{k=0}^{m}$}(a-(b-1-k)i)(a+(b-k-2)i)a_kv_{0,k}.
\end{equation}
Since we can choose suitable $i,k\in\Z$ such that
$$(a-(b-1-k)i)(a+(b-k-2)i)a_k\not=0,$$ then if we apply $\rho^0=id, \rho, \rho^2,\cdots,
\rho^m$ to $v$ respectively, we obtain $m+1$ equations concerning
$v_{0,k}\,(0\leq k\leq m)$ correspondingly. Denote the matrix of the
coefficients of $v_{0,k}\,(0\leq k\leq m)$ by $M$ and let
$\Delta'=\det\,M.$ By linear algebra, it is easy to see that
$\Delta'\not=0.$ Since by our assumption, $V'$ is a submodule of
$A_{a,b}$, we know that $\rho^{n}v\in V'$ for $0\leq n\leq m,$ then
it follows that $v_{0,k}\in V'$ for all $0\leq k\leq m$. Therefore
by index shifting, we can suppose $v_{0,0}\in V'$. Our aim is to
find out under what condition we have $V'=A_{a,b}$.

At first, we consider the case $a\notin\Z$. By (\ref{A-a-b}), we
know that
\begin{equation}\label{irred1}
v_{i,j}\in V'  {\rm \ \ for} \ i,j\in\Z {\rm \ and}\
(j+1)a+(b-1)i\not=0.
\end{equation}
So in this case we only need to consider whether $v_{i,j}\in V'$ for
$(j+1)a+(b-1)i=0$. Since $a\notin\Z$, by (\ref{irred1}), we get that
$v_{0,1}\in V'$. Applying $L_{i,j-1}$ to $v_{0,1}$, using
(\ref{A-a-b}), we obtain that
\begin{eqnarray}\label{irred2}
\!\!\!\!\!\!\!\!\!\!\!\!&\!\!\!\!\!\!\!\!\!\!\!\!\!\!\!&
L_{i,j-1}v_{0,1}=(j a+(b-2)i)v_{i,j}=-(a+i)v_{i,j},
\end{eqnarray}
for $i,j\in\Z$ and $(j+1)a+(b-1)i=0$. Since $a\notin\Z$,
$a+i\not=0$, we have that $v_{i,j}\in V'$ for $(j+1)a+(b-1)i=0$.
This together with (\ref{irred1}) shows that $V'=A_{a,b}$.
Therefore, if $a\not\in\Z$, then  $A_{a,b}$ is irreducible.

Now we consider the case $a\in\Z$ and $b\notin\Z$. By (\ref{A-a-b}),
we have $v_{i,j}\in V'$ for $i,j\in\Z$ and $i\not=0$. Applying
$L_{-i,0}$ to $v_{i,j}$, using (\ref{A-a-b}), we obtain that
\begin{equation}\label{irred3}
L_{-i,0}v_{i,j}=(a+i-(b-j-2)i)v_{0,j}  {\rm \ \ for} \ i,j\in\Z {\rm
\ and}\ i\not=0.
\end{equation}
Since $a\in\Z$ and $b\notin\Z$, it follows that $v_{0,j}\in V'$.
Thus we have $V'=A_{a,b}$. Therefore, if $a\in\Z$ and $b\not\in\Z$,
then $A_{a,b}$ is irreducible. Obviously, there is no deformation
when $A_{a,b}$ is irreducible.

Finally we consider the case $a,b\in\Z$. By (\ref{A-a-b}), we have
\begin{equation}\label{reducible}
L_{i,j}v_{-a,b-1}=0  {\rm \ \ for\ any} \ i,j\in\Z {\rm \ and}\
a,b\in\Z.
\end{equation}
Set $V_1={\rm span}{\{v_{-a,b-1}\}}$, $V_2={\rm
span}{\{v_{-a,b-2}\}}$, $V_3={\rm
span}{\{v_{k,l}\,|\,(k,l)\not=(-a,b-1),(-a,b-2)\}}$. By
(\ref{A-a-b}), we know that $A'_{a,b}:=V_1$ is a 1-dimensional
$\LL-$submodule of $A_{a,b}$. Since $A_{a,b}$ is indecomposable,
then by (\ref{A-a-b}), we obtain that $A_{a,b}$ has three
composition factors, the first is the 1-dimensional submodule
$A'_{a,b}$, the second is the quotient module $A''_{a,b}:=V_3/V_1$
and the third is the quotient module $A'''_{a,b}:=V_2/V_3$.
Therefore (i) and (ii) of Theorem $\ref{theo}(2)$ follows.

Obviously, when $a,b\in\Z$, we have $A_{a,b}\cong A_{0,0}$ by index
shifting. Since $A_{0,0}$ is reducible, then it might have
deformations. In order to find out all deformations, we need to find
out all possible different actions. For $i,j,k,l\in\Z$, suppose
\begin{eqnarray}\label{general}
\!\!\!\!\!\!\!\!\!\!\!\!&\!\!\!\!\!\!\!\!\!\!\!\!\!\!\!&
L_{i,j}v_{k,l}=((j+1)k-(l+1)i)v_{i+k,j+l} {\rm \ \ \ \ for} \
(k,l),(i+k,j+l)\not=(0,-1),(0,-2),\nonumber
\\\!\!\!\!\!\!\!\!\!\!\!\!&\!\!\!\!\!\!\!\!\!\!\!\!\!\!\!&
L_{i,j}v_{0,-1}=a_{i,j}v_{i,j-1} {\rm \ \ \ \ for\ some}\
a_{i,j}\in\F,\nonumber
\\\!\!\!\!\!\!\!\!\!\!\!\!&\!\!\!\!\!\!\!\!\!\!\!\!\!\!\!&
L_{i,j}v_{0,-2}=b_{i,j}v_{i,j-2} {\rm \ \ \ \ for\ some}\
b_{i,j}\in\F,
\\\!\!\!\!\!\!\!\!\!\!\!\!&\!\!\!\!\!\!\!\!\!\!\!\!\!\!\!&
L_{i,j}v_{-i,-j-1}=c_{i,j}v_{0,-1} {\rm \ \ \ \ for\ some}\
c_{i,j}\in\F,\nonumber
\\\!\!\!\!\!\!\!\!\!\!\!\!&\!\!\!\!\!\!\!\!\!\!\!\!\!\!\!&
L_{i,j}v_{-i,-j-2}=d_{i,j}v_{0,-2} {\rm \ \ \ \ for\ some}\
d_{i,j}\in\F.\nonumber
\end{eqnarray}
For convenience, in the following, we also denote $V_1={\rm
span}{\{v_{0,-1}\}}$, $V_2={\rm span}{\{v_{0,-2}\}}$ and $V_3={\rm
span}{\{v_{k,l}\,|\,(k,l)\not=(0,-1),(0,-2)\}}$.

Now we only need to consider the following cases.

{\it Case 1: $v_{0,-2}=0$ and $V_1$ is a submodule of $V$.}

In this case, for $i,j,k,l\in\Z$, we can suppose
\begin{eqnarray}\label{suppose3}
\!\!\!\!\!\!\!\!\!\!\!\!&\!\!\!\!\!\!\!\!\!\!\!\!\!\!\!&
L_{i,j}v_{k,l}=((j+1)k-(l+1)i)v_{i+k,j+l} {\rm \ \ for} \
(k,l),(i+k,j+l)\not=(0,-1),
\\\!\!\!\!\!\!\!\!\!\!\!\!&\!\!\!\!\!\!\!\!\!\!\!\!\!\!\!&
L_{i,j}v_{0,-1}=0 {\rm \ \ for}\ (i,j)\not=(0,0),\label{suppose4'}
\\\!\!\!\!\!\!\!\!\!\!\!\!&\!\!\!\!\!\!\!\!\!\!\!\!\!\!\!&
L_{i,j}v_{-i,-j-1}=c_{i,j}v_{0,-1} {\rm \ \ for\ some}\
c_{i,j}\in\F.\label{suppose4}
\end{eqnarray}
Then we have
\begin{eqnarray}\label{triple}
\!\!\!\!\!\!\!\!\!\!\!\!&\!\!\!\!\!\!\!\!\!\!\!\!\!\!\!&
L_{i,j}L_{1,0}v_{-i-1,-j-1}=(j-i-1)c_{i,j}v_{0,-1}  {\rm \ \ for} \
(i,j)\not=(0,0),
\nonumber\\\!\!\!\!\!\!\!\!\!\!\!\!&\!\!\!\!\!\!\!\!\!\!\!\!\!\!\!&
L_{i,j}L_{0,1}v_{-i,-2-j}=-2ic_{i,j}v_{0,-1}  {\rm \ \ for} \
(i,j)\not=(0,0),
\nonumber\\\!\!\!\!\!\!\!\!\!\!\!\!&\!\!\!\!\!\!\!\!\!\!\!\!\!\!\!&
L_{i,j}L_{-i,0}v_{0,-j-1}=-ijc_{i,j}v_{0,-1}  {\rm \ \ for} \
(i,j)\not=(0,0).
\end{eqnarray}
Applying $L_{i,j}L_{1,0}=(j+1-i)L_{i+1,j}+L_{1,0}L_{i,j}$,
$L_{i,j}L_{0,1}=-2iL_{i,j+1}+L_{0,1}L_{i,j}$ and
$L_{i,j}L_{-i,0}=-i(j+2)L_{0,j}+L_{-i,0}L_{i,j}$  to
$v_{-i-1,-j-1}$, $v_{-i,-j-2}$ and $v_{0,-j-1}$ respectively,
comparing the coefficients of $v_{0,-1}$, using (\ref{suppose3}),
(\ref{suppose4}), (\ref{suppose4'}) and (\ref{triple}), we can
deduce the following formulas:
\begin{eqnarray}\label{triple-1}
\!\!\!\!\!\!\!\!\!\!\!\!&\!\!\!\!\!\!\!\!\!\!\!\!\!\!\!&
(j+1-i)c_{i+1,j}-(i+j+1)c_{1,0}=(j-i-1)c_{i,j}  {\rm \ \ for} \
(i,j)\not=(0,0),
\\\!\!\!\!\!\!\!\!\!\!\!\!&\!\!\!\!\!\!\!\!\!\!\!\!\!\!\!&
c_{0,0}=0, \ \ \ c_{i,j}=c_{i,0}  {\rm \ \ for} \ i,j\in\Z {\rm \
and\ }i\not=0,
\label{triple-2}\\\!\!\!\!\!\!\!\!\!\!\!\!&\!\!\!\!\!\!\!\!\!\!\!\!\!\!\!&
(j+2)c_{0,j}-jc_{-i,0}=jc_{i,j}  {\rm \ \ for} \ i,j\in\Z {\rm \
and\ }i\not=0.\label{triple-3}
\end{eqnarray}
Setting $i=0$ and $j\not=0$ in (\ref{triple-1}), we have
\begin{equation}\label{bee}
(j+1)c_{1,j}-(j+1)c_{1,0}=(j-1)c_{0,j}.
\end{equation}
Using (\ref{triple-2}), we know that $c_{1,j}=c_{1,0}$ for any
$j\in\Z$. Then (\ref{bee}) becomes $(j-1)c_{0,j}=0$ for $j\in\Z$ and
$j\not=0$. Therefore, by (\ref{triple-2}), we can get
\begin{equation}\label{help}
c_{0,j}=0 {\rm \ \ for\ } j\in\Z {\rm \ and}\ j\not=1.
\end{equation}
Taking $j=2$ in (\ref{triple-3}), using (\ref{triple-2}) and
(\ref{help}), we can deduce that
\begin{equation}\label{back}
c_{-i,0}=-c_{i,0} {\rm \ \ for\ all\ } i\in\Z.
\end{equation}
Setting $j=1$ in (\ref{triple-3}), we have
$3c_{0,1}-c_{-i,0}=c_{i,0}$. By (\ref{back}), it follows that
$c_{0,1}=0$. This together with (\ref{help}) show that $c_{0,j}=0$
for all $j\in\Z$. Then (\ref{triple-2}) becomes
\begin{equation}\label{fact-1}
c_{i,j}=c_{i,0} {\rm \ \ for\ any\ } i,j\in\Z.
\end{equation}
And (\ref{triple-1}) becomes
\begin{equation}\label{fact-2}
(j+1-i)c_{i+1,0}-(i+j+1)c_{1,0}=(j-i-1)c_{i,0} {\rm \ \ for\ any\ }
i,j\in\Z.
\end{equation}
Setting $j=i-1$ in (\ref{fact-2}), it gives that $c_{i,0}=ic_{1,0}$
for any $i\in\Z$. This together with (\ref{fact-1}) shows that
$c_{i,j}=ic_{1,0}.$ If $c_{1,0}\ne0$, by rescaling $v_{0,0}$, we can
suppose $c_{1,0}=-1$. Therefore in this case, we obtain that
$c_{i,j}=-i$.

{\it Case 2: $v_{0,-1}=0$ and $V_2$ is a submodule of $V$.}

In this case, for $i,j,k,l\in\Z$, we can suppose
\begin{eqnarray}\label{case2}
\!\!\!\!\!\!\!\!\!\!\!\!&\!\!\!\!\!\!\!\!\!\!\!\!\!\!\!&
L_{i,j}v_{k,l}=((j+1)k-(l+1)i)v_{i+k,j+l} {\rm \ \ for} \
(k,l),(i+k,j+l)\not=(0,-2),
\\\!\!\!\!\!\!\!\!\!\!\!\!&\!\!\!\!\!\!\!\!\!\!\!\!\!\!\!&
L_{i,j}v_{0,-2}=0 {\rm \ \ for}\ (i,j)\not=(0,0),
\\\!\!\!\!\!\!\!\!\!\!\!\!&\!\!\!\!\!\!\!\!\!\!\!\!\!\!\!&
L_{i,j}v_{-i,-j-2}=d_{i,j}v_{0,-2} {\rm \ \ for\ some}\
d_{i,j}\in\F.
\end{eqnarray}

Using the same methods developed in Case 1, we can deduce that
$d_{i,j}=i.$

{\it Case 3: $v_{0,-2}=0$ and $V_3$ is a submodule of $V$.}

For $i,j,k,l\in\Z$, in this case we can suppose
\begin{eqnarray}\label{suppose1}
\!\!\!\!\!\!\!\!\!\!\!\!&\!\!\!\!\!\!\!\!\!\!\!\!\!\!\!&
L_{i,j}v_{k,l}=((j+1)k-(l+1)i)v_{i+k,j+l} {\rm \ \ for} \ (k,l),
(i+k,j+l)\not=(0,-1),
\\\!\!\!\!\!\!\!\!\!\!\!\!&\!\!\!\!\!\!\!\!\!\!\!\!\!\!\!&
L_{i,j}v_{0,-1}=a_{i,j}v_{i,j-1} {\rm \ \ for\ some}\
a_{i,j}\in\F,\label{suppose2}
\\\!\!\!\!\!\!\!\!\!\!\!\!&\!\!\!\!\!\!\!\!\!\!\!\!\!\!\!&
L_{i,j}v_{-i,-j-1}=0 {\rm \ \ for}\
(i,j)\not=(0,0).\label{suppose2'}
\end{eqnarray}
Then we have
\begin{eqnarray}\label{triple'}
\!\!\!\!\!\!\!\!\!\!\!\!&\!\!\!\!\!\!\!\!\!\!\!\!\!\!\!&
L_{1,0}L_{i,j}v_{0,-1}=(i-j)a_{i,j}v_{i+1,j-1}  {\rm \ \ for} \
(i,j)\not=(0,0),(-1,0),
\nonumber\\\!\!\!\!\!\!\!\!\!\!\!\!&\!\!\!\!\!\!\!\!\!\!\!\!\!\!\!&
L_{-1,0}L_{i,j}v_{0,-1}=(i+j)a_{i,j}v_{i-1,j-1} {\rm \ \ for} \
(i,j)\not=(0,0),(1,0),
\nonumber\\\!\!\!\!\!\!\!\!\!\!\!\!&\!\!\!\!\!\!\!\!\!\!\!\!\!\!\!&
L_{0,1}L_{i,j}v_{0,-1}=2ia_{i,j}v_{i,j}  {\rm \ \ for} \
(i,j)\not=(0,0),(0,-1).
\end{eqnarray}
Applying $L_{i,j}L_{1,0}=(j+1-i)L_{i+1,j}+L_{1,0}L_{i,j}$,
$L_{i,j}L_{0,1}=-2iL_{i,j+1}+L_{0,1}L_{i,j}$ and
$L_{i,j}L_{-1,0}=-(i+j+1)L_{i-1,j}+L_{-1,0}L_{i,j}$  to $v_{0,-1}$,
comparing the coefficients of $v_{i+1,j-1}$, $v_{i,j}$ and
$v_{i-1,j-1}$ respectively, using (\ref{suppose1}),
(\ref{suppose2}), (\ref{suppose2'}) and (\ref{triple'}), we can
deduce the following formulas:
\begin{eqnarray}\label{triple'-1}
\!\!\!\!\!\!\!\!\!\!\!\!&\!\!\!\!\!\!\!\!\!\!\!\!\!\!\!&
(j-i+1)a_{i+1,j}-(j-i)a_{i,j}=(j+1)a_{1,0} {\rm \ \ for} \
(i,j)\not=(-1,0),
\\\!\!\!\!\!\!\!\!\!\!\!\!&\!\!\!\!\!\!\!\!\!\!\!\!\!\!\!&
a_{0,0}=0, \ \ \ \ a_{i,j}=a_{i,0}+\frac{j}{2}a_{0,1}  {\rm \ \ for}
\ i,j\in\Z {\rm \ and\ }i\not=0,
\label{triple'-3}\\\!\!\!\!\!\!\!\!\!\!\!\!&\!\!\!\!\!\!\!\!\!\!\!\!\!\!\!&
(i+j)a_{i,j}-(i+j+1)a_{i-1,j}=-(j+1)a_{-1,0}  {\rm \ \ for} \
(i,j)\not=(1,0).\label{triple'-4}
\end{eqnarray}
Replacing $i$ by $i-1$ in (\ref{triple'-1}), we can get
$(j-i+2)a_{i,j}-(j-i+1)a_{i-1,j}=(j+1)a_{1,0}  {\rm \ for} \
(i,j)\not=(-2,0).$ From this equation and (\ref{triple'-3}),
(\ref{triple'-4}), we obtain that
\begin{equation}\label{main-1}
(j+1)(2a_{i,j}+(i-j-1)a_{-1,0}-(i+j+1)a_{1,0})=0 {\rm \ \ for\ }
(i,j)\not=(-2,0).
\end{equation}
Setting $(i,j)=(0,0)$ in (\ref{main-1}), by (\ref{triple'-3}), we
have $a_{-1,0}=-a_{1,0}.$ Then (\ref{main-1}) becomes
\begin{equation}\label{main-2}
(j+1)(2a_{i,j}-2ia_{1,0})=0 {\rm \ \ for\ } (i,j)\not=(-2,0).
\end{equation}
Letting $(i,j)=(0,1)$ in (\ref{main-2}), it gives that $a_{0,1}=0$.
Then (\ref{triple'-3}) becomes
\begin{equation}\label{main-3}
a_{i,j}=a_{i,0} {\rm \ \ for\ } i\not=0.
\end{equation}
The above formula together with (\ref{main-2}) shows that
$a_{i,j}=ia_{1,0}$ for all $i,j\in\Z$. Therefore, as before, in this
case we can obtain that $a_{i,j}=i$.

{\it Case 4: $v_{0,-1}=0$ and $V_3$ is a submodule of $V$.}

For $i,j,k,l\in\Z$, in this case we can suppose
\begin{eqnarray}\label{suppose1}
\!\!\!\!\!\!\!\!\!\!\!\!&\!\!\!\!\!\!\!\!\!\!\!\!\!\!\!&
L_{i,j}v_{k,l}=((j+1)k-(l+1)i)v_{i+k,j+l} {\rm \ \ for} \ (k,l),
(i+k,j+l)\not=(0,-2),
\\\!\!\!\!\!\!\!\!\!\!\!\!&\!\!\!\!\!\!\!\!\!\!\!\!\!\!\!&
L_{i,j}v_{0,-2}=b_{i,j}v_{i,j-2} {\rm \ \ for\ some}\
b_{i,j}\in\F,\label{suppose2}
\\\!\!\!\!\!\!\!\!\!\!\!\!&\!\!\!\!\!\!\!\!\!\!\!\!\!\!\!&
L_{i,j}v_{-i,-j-2}=0 {\rm \ \ for}\
(i,j)\not=(0,0).\label{suppose2'}
\end{eqnarray}

Using the same methods as in Case 3, we can deduce that $b_{i,j}=i$.

{\it Case 5: $v_{0,-1}\not=0$, $v_{0,-2}\not=0$ and $V_1$, $V_2$ are
submodules of $V$.}

From the discussions in Case 1 and Case 2, in this case we can
obtain a deformation of the adjoint module $A_{0,0}$ of $\LL$,
denoted by $A$, and so (\ref{B-a-b}) in Theorem \ref{theo} holds (it
is straightforward to verify that $A$ is a module).

{\it Case 6: $v_{0,-1}\not=0$, $v_{0,-2}\not=0$ and $V_3$ is a
submodule of $V$.}

From the discussions in Case 3 and Case 4, in this case we can
obtain a deformation of the adjoint module $A_{0,0}$ of $\LL$,
denoted by $B$, and so (\ref{B-a-b-c}) in Theorem \ref{theo} holds.

{\it Case 7: $v_{0,-1}\not=0$, $v_{0,-2}\not=0$ and $V_2$ is a
submodule of $V$.}

From the discussions in Case 2 and Case 3, as before, in this case
we can obtain a deformation of the adjoint module $A_{0,0}$ of
$\LL$, denoted by $C$, and so (\ref{B-a-b-d}) in Theorem \ref{theo}
holds.

From (\ref{B-a-b}), (\ref{B-a-b-c}) and (\ref{B-a-b-d}), it is easy
to see that $\rm(iii)$, $\rm(iv)$ and $\rm(v)$ of Theorem
\ref{theo}$(2)$ hold. This completes the proof of Theorem
\ref{theo}.\hfill$\Box$\vskip5pt

\cl{\bf\small Acknowledgement}

 \par\

 {\small The author would like
to thank Dr. Rencai Lu for pointing out some errors in the first
version. }

\vskip12pt

\cl{\bf References}\vs{0pt}

\vskip5pt\small
\parindent=8ex\parskip=1pt\baselineskip=1pt\parskip0.001in

\re{DZ} D. Dokovic, K. Zhao, Derivations, isomorphisms and
cohomology of generalized Block algebras, {\it Algebra Colloq.} {\bf
3} (1996), 245--272.

\re{LT}  W. Lin, S.  Tan, Nonzero level Harish-Chandra modules over
the Virasoro-like algebra, {\it J. Pure Appl. Algebra} {\bf 204}
(2006), 90--105.

\re{S1} Y. Su, Quasifinite representations of a Lie algebra of Block
type, {\it J. Algebra} {\bf 276} (2004), 117--128.

\re{S2} Y. Su, Quasifinite representations of a family of Lie
algebras of Block type, {\it J. Pure Appl. Algebra} {\bf 192}
(2004), 293--305.

\re{S3} Y. Su, Simple modules over the high rank Virasoro algebras,
{\it Comm.  Algebra} {\bf 29} (2001), 2067--2080.

\re{S4} Y. Su,  Classification of Harish-Chandra modules over the
higher rank Virasoro algebras, {\it Comm. Math. Phys.} {\bf240}
(2003), 539--551.

\re{WZ} X. Wang, K. Zhao, Verma modules over the Virasoro-like
algebra, {\it J. Australia Math.} {\bf80} (2006), 179--191.

\re{WS} Y. Wu, Y. Su, Highest weight representations of a Lie
algebra of Block type, to appear in {\it Science in China A}
(math.QA/0511733).

\re{X1} X. Xu, Generalizations of Block algebras, {\it Manuscripta
Math}. {\bf 100} (1999), 489--518.

\re{X2} X. Xu, Quadratic conformal superalgebras, {\it J. Algebra}
 {\bf 231} (2000), 1--38.

\re{ZZ} H. Zhang, K. Zhao, Representations of the Virasoro-like Lie
algebra and its q-analog, {\it Comm. Algebra} {\bf 24} (1996),
4361--4372.

\re{Z} K. Zhao, Weight modules over generalized Witt algebras with
 1-dimensional weight spaces, {\it Forum Math.} {\bf16} (2004),
 725--748.

\re{ZM} L. Zhu, D. Meng, Structure of degenerate Block algebras,
{\it Algebra Colloq.} {\bf 10} (2003), 53--62.

\end{document}